\documentclass{amsart}

\usepackage{amssymb}

\newtheorem{thm}{Theorem}%[section]

\newtheorem{lem}[thm]{Lemma}
\newtheorem{cor}[thm]{Corollary}

\newtheorem{prop}[thm]{Proposition}

\newtheorem{principle}[thm]{Principle} %!!!!!!!!!!!!!!!!!!!!!!
\theoremstyle{definition}
\newtheorem{defn}[thm]{Definition}

\newtheorem{say}[thm]{}
\newtheorem{exmp}[thm]{Example}

\newtheorem{ques}[thm]{Question}    %!!!!!!!!!!!!!!!!!!!!

\newtheorem{rem}[thm]{Remark}          

\newtheorem{ack}{Acknowledgments}

\newtheorem{defn-thm}[thm]{Definition--Theorem}  %!!!!!!!!!!!!!!!!!!!!!!!!
\newtheorem{defn-lem}[thm]{Definition--Lemma}  %!!!!!!!!!!!!!!!!!!!!!!!!
\theoremstyle{remark}
\newtheorem{claim}[thm]{Claim}

\renewcommand{\c}[0]{{\mathbb C}}  

\renewcommand{\o}[0]{{\mathcal O}} 
\newcommand{\z}[0]{{\mathbb Z}}

\newcommand{\h}[0]{{\mathbb H}}
\newcommand{\p}[0]{{\mathbb P}}

\newcommand{\map}[0]{\dasharrow}
\newcommand{\qtq}[1]{\quad\mbox{#1}\quad}
\newcommand{\spec}[0]{\operatorname{Spec}}
\newcommand{\pic}[0]{\operatorname{Pic}}

\newcommand{\rank}[0]{\operatorname{rank}}

\newcommand{\supp}[0]{\operatorname{Supp}}

\newcommand{\im}[0]{\operatorname{im}}

\newcommand{\Hom}[0]{\operatorname{Hom}}

\newcommand{\sing}[0]{\operatorname{Sing}}

\newcommand{\onto}[0]{\twoheadrightarrow}

\newcommand{\diag}[0]{\operatorname{diag}}

\newcommand{\pardeg}[0]{\operatorname{pardeg}}
\newcommand{\tD}[0]{{\tilde{D}}}
\newcommand{\weight}[0]{\operatorname{weight}}

\def\into{\DOTSB\lhook\joinrel\to}

\begin{document}
\bibliographystyle{amsalpha}

\title{Holonomy groups of stable vector bundles}
\author{V. Balaji and J\'anos Koll\'ar}

\today
\maketitle

%\tableofcontents

Let $M$ be a Riemannian manifold and $E$ a vector bundle
with a connection  $\nabla$. Parallel transport
along loops gives a representation
of  the loop group of $M$ with base point $x$
into the orthogonal group  $O(E_x)$ of the fiber at $x$
(see, for instance,  \cite{k-n, bryant}).

If $X$ is a complex manifold and $E$ a holomorphic vector bundle,
then usually there are no holomorphic connections on $E$.
One can, nonetheless,  define a close analog of the
holonomy representation in the complex setting if
$E$ is a stable vector bundle and $X$ is projective algebraic.

Assume for simplicity that $c_1(E)=0$,
a condition that we remove later.
By  Mehta--Ramanathan  \cite{meh-ram1},
if $x\in C\subset X$ is a sufficiently general
complex curve, then $E|_C$ is  also stable
and so by  a result  of
Narasimhan--Seshadri \cite{nar-ses}
it corresponds to a unique unitary representation
$\rho:\pi_1(C)\to U(E_x)$.
We call it the Narasimhan--Seshadri representation  of $E|_C$.

The image of the representation, and even the Hermitian
form on $E_x$ implicit in its definition, depend on the
choice of $C$, but the picture stabilizes
if we look at the Zariski closure 
of the image in $GL(E_x)$.
The resulting group can also be characterized in  different ways.

\begin{thm} \label{hol.thm.intr}
Let $X$ be a smooth projective variety,
$H$ an ample divisor on $X$, $E$ a stable
vector bundle with $\det E\cong \o_X$ and $x\in X$ a point.
Then there is a unique reductive subgroup
$H_x(E)\subset SL(E_x)$, called the {\em holonomy} group of $E$,
characterized by either of  the  two properties:
\begin{enumerate}
\item $H_x(E)\subset SL(E_x)$ is the smallest algebraic subgroup
satisfying the following:

For every curve  $x\in C\subset X$ such that $E|_C$ is stable,
the image of the Narasimhan--Seshadri representation
$\rho:\pi_1(C)\to U(E_x)$ is contained in $H_x(E)$.
\item If $C$ is  sufficiently general, then
the image of the Narasimhan--Seshadri representation
is Zariski dense in $H_x(E)$.
\end{enumerate}
Furthermore:
\begin{enumerate}\setcounter{enumi}{2}
\item For every $m,n$, the
fiber map $F\to F_x$  gives 
 a one--to--one correspondence between
direct summands of $E^{\otimes m}\otimes (E^*)^{\otimes n} $ and 
$H_x(E)$-invariant subspaces of $E_x^{\otimes m}\otimes (E_x^*)^{\otimes n}$.
\item  The conjugacy class of $H_x(E)$ is the smallest  
reductive   conjugacy class   $G$
such that the structure group of $E$ can be reduced to $G$.
\end{enumerate}
\end{thm}

\begin{rem}\label{main.thm.rem}
  (1) The existence of a smallest reductive structure group
is established in \cite[Thm.2.1]{bog-sb}.

(2) We emphasize that the holonomy group is defined
as a subgroup of $GL(E_x)$ and not just as a conjugacy class
of subgroups.

(3) It follows from  (\ref{hol.thm.intr}.1) 
 that the holonomy group does not
depend on $H$. Thus the definition of the holonomy group
makes sense for any vector bundle that is stable
with respect to some ample divisor $H$.

(4) The property (\ref{hol.thm.intr}.3) almost characterizes
the holonomy group. The only remaining ambiguity
comes from the center of $GL(E_x)$. In general,
the holonomy group is determined by knowing, for every $m,n\geq 0$, the
direct summands of 
$E^{\otimes m}\otimes (E^*)^{\otimes n} $ 
and also knowing which rank 1 summands are isomorphic to $\o_X$.

(5) The above theorem has immediate generalizations to the
 case when $X$ is a normal variety, $E$ a reflexive sheaf
with arbitrary $\det E$ or a  sheaf with parabolic structure.
 These are discussed in 
 (\ref{hol.thm}) and  (\ref{parab.holom.sefd}).
The case of Higgs bundles will be considered elsewhere.

(6) For some closely related ideas and applications to
the construction of stable principal bundles on surfaces, see
\cite{balaji}.
\end{rem}

Our next aim is to study and use holonomy groups
by relying on the following:

\begin{principle}\label{hsm.princ} 
Let $E$ be a stable vector bundle on a smooth projective variety
$X$.
\begin{enumerate} 
\item If $E$ is ``general'' then
 the holonomy group $H_x(E)$ is ``large'',
meaning, for instance, that   $H_x(E)\supset SL(E_x)$.
\item Otherwise 
 there is geometric reason why $H_x(E)$ is small.
\end{enumerate}
\end{principle}

Let  $\rho:\pi_1(X)\to U(V)$ be an irreducible
 representation with finite image
$G$ and  $E_{\rho}$  the corresponding flat
vector  bundle on $X$.
Then $H_x(E_{\rho})=G$. Understanding
$G$ in terms of its representations $V^{\otimes m}$ is certainly
possible, but  it quickly leads to  intricate
questions of finite group theory.
(See \cite{g-t} for such an example.)
There is a significant case
when we can avoid the complications coming from finite
subgroups of $GL(V)$.

\begin{prop} \label{hol.simpconn.prop}
If $X$ is simply connected
then $H_x(E)$ is connected.
\end{prop}

The representation theory of connected reductive groups
is quite well understood, and this enables us to get some
illustration of the above Principle.

\begin{prop}\label{S^m.cond}
Let $E$ be a stable vector bundle on a
 simply connected smooth projective variety $X$.
Then the following are equivalent:
\begin{enumerate}
\item $S^mE$ is stable (that is, indecomposable) for some $m\geq 2$.
\item $S^mE$ is stable (that is, indecomposable) for every $m\geq 2$.
\item The holonomy is one of the following: 
\begin{enumerate}
\item $SL(E_x)$ or $GL(E_x)$, 
\item $Sp(E_x)$ or $GSp(E_x)$ 
for a suitable nondegenerate  symplectic form 
on $E_x$ (and  $\rank E$ is even).
\end{enumerate}
\end{enumerate}
\end{prop}

Note that the statements
(\ref{S^m.cond}.1) and (\ref{S^m.cond}.2)
do not involve the holonomy group, but it
 is not clear to us how to prove their equivalence
without using holonomy.

If $X$ is not simply connected, 
the results of  \cite{g-t2} imply the
following:

\begin{cor}\label{gural.cor} Let $E$ be a stable vector bundle on a
smooth projective variety $X$ of rank $\geq 7$. Then the following are
equivalent:
\begin{enumerate}
\item $S^rE$ is stable  for some $r\geq 4$.
\item $S^rE$ is stable  for every $r\geq 2$.
\item The commutator of the  holonomy group is either $SL(E_x)$ or $Sp(E_x)$. 
\end{enumerate}
\end{cor}

The  rank $\leq 6$ cases  can, in principle,  be enumerated by hand.
As Guralnick communicated to us, there are probably only two examples where 
(\ref{gural.cor}) fails. One is the icosahedral group
in $SL(2,\c)$ and the other is a double cover of the second Janko group
$2J_2$ acting on $\c^6$.

Another illustration of the Principle (\ref{hsm.princ})
is the following partial description of  
low rank bundles.

\begin{prop}\label{low.rank.prop}  Let $E$ be a stable vector bundle on a
 simply connected smooth projective variety $X$.
Assume that $\det E\cong \o_X$ 
and  $\rank E\leq 7$.
Then one of the following holds.
\begin{enumerate}
\item The holonomy group is $SL(E_x)$.
\item The holonomy is contained in $SO(E_x)$ or $Sp(E_x)$.
In particular, $E\cong E^*$ and the odd Chern classes of $E$ are 2-torsion.
\item $E$ is obtained from a rank 2 vector bundle $F_2$
and a rank 3 vector bundle $F_3$.
There are 2 such cases which are neither orthogonal nor
symplectic:
\begin{enumerate}
\item $\rank E=6$ and $E\cong S^2F_3$, or
 \item $\rank E=6$ and $E\cong F_2\otimes F_3$.
\end{enumerate}
\end{enumerate}
\end{prop}

There are two reasons why 
a result of this type gets more complicated
for higher rank bundles.

First,  already in rank 7, we have vector bundles with $G_2$ holonomy. 
It is not on our list separately since $G_2\subset SO_7$.
It is quite likely that there is some very nice geometry
associated  with $G_2$ holonomy, but this remains to be
discovered.
Similarly, the other exceptional groups must all appear
for the higher rank cases.

Second, and this is more serious,
there are many cases where the holonomy group is
not simply connected, for instance
 $PGL$. In this case there is a Brauer
obstruction to lift the structure group to $GL$
and to write $E$ in terms of a lower rank bundle using
representation theory.
We study this in (\ref{Gl-reps.brauer}).
In the low rank cases we are saved by the accident that
such representations happen to be either orthogonal or
symplectic, but this definitely fails in general.

\begin{say}[Comparison with the differential geometric holonomy]

For  the tangent bundle of a smooth projective variety $X$,
one gets two notions of holonomy.
The classical differential geometric holonomy
and the algebraic holonomy defined earlier.
These are related in some ways, but the precise
relationship is still unclear.

First of all, the algebraic holonomy makes sense whenever
$T_X$ is stable, and it does not depend on the choice of
a metric on $X$. The differential geometric holonomy
depends on the metric chosen.

If $X$ admits a  K\"ahler--Einstein metric, then
 its holonomy group, which
 is a subgroup of the unitary group $U(T_xX)$,
is canonically associated to $X$.

By contrast, the algebraic holonomy is not unitary.
For a general curve $C\subset X$, the
 Narasimhan--Seshadri representation
gives a subgroup of a unitary group, but the
Hermitian form defining the unitary group in question does depend on  $C$,
except when $X$ is a quotient of an Abelian variety.

Thus the processes that define the holonomy group
in algebraic geometry and in differential geometry
are quite different.
It is, nonetheless, possible, that the two 
holonomy groups are closely related.

\begin{ques} \label{hol=hol?} 
Let $X$ be a simply connected smooth projective variety
which has a K\"ahler--Einstein metric. Is the
algebraic holonomy group of $T_X$  the  complexification
of the  differential geometric holonomy group?
\end{ques}

For non simply connected varieties 
the  differential geometric holonomy group may have infinitely many
connected components, and
one may need to take
the complexification of its Zariski closure instead.

It is possible that (\ref{hol=hol?}) holds for the simple reason that 
the algebraic holonomy group of a tangent bundle is almost
always $GL_n$. 
The differential geometric holonomy group is
almost always $U_n$, with two notable exceptions. In both
of them, the answer to (\ref{hol=hol?}) is positive.

\begin{prop}[Calabi--Yau varieties]\label{hol=hol.CY} Let $X$ be a 
simply connected smooth projective variety $X$
such that $K_X=0$ which is not a direct product.

The differential geometric holonomy group is 
either $SU_n(\c)$ or $U_n(\h)$. Correspondingly, the algebraic holonomy group
is $SL_n(\c)$ (resp.\ $Sp_{2n}(\c)$).
\end{prop}

\begin{prop}[Homogeneous spaces] \label{G/P.prop} 
Let $X=G/P$ be a homogeneous space such that
the stabilizer representation of $P$ on $T_xX$ is irreducible.
Then $T_X$ is stable and the algebraic holonomy
is the image of this stabilizer representation.
\end{prop}

\end{say}

\section{Variation of monodromy groups}

\begin{say} Let $C$ be a smooth projective curve over $\c$
and $x\in C$ a point. Every unitary representation
$\rho: \pi_1(C, x)\to U(\c^r)$ gives a flat vector bundle $E_{\rho}$
of rank $r$. By \cite{nar-ses}, this gives a
real analytic one--to--one correspondence
between conjugacy classes of  unitary representations
and polystable vector bundles of rank $r$ and degree $0$.

The similar correspondence between representations
and polystable vector bundles of rank $r$ and degree $d\neq 0$
is less natural and it depends on an additional point of $C$.

 Let $C$ be a smooth projective curve over $\c$
and $x\neq c\in C$ two points.
Let $\Gamma\subset \pi_1(C\setminus c, x)$ denote the conjugacy class
consisting of counterclockwise lassos around $c$.

A unitary representation
$$
\rho: \pi_1(C\setminus c, x)\to U(\c^r) 
\qtq{such that}\rho(\gamma)=e^{2\pi i d/r}{\mathbf 1}
$$ 
 for every $\gamma\in \Gamma$
is said to have type $d/r$.
(In the original definition this is called type $d$.
Using type $d/r$ has the advantage that
irreducible subrepresentations have the same type.)
Note that the type is well defined only modulo $1$.

By \cite{nar-ses}, for every $r$ and $d$ the following hold: 
\begin{enumerate}
\item There is a one--to-one correspondence 
$$
NS:(C,c,x,E)\mapsto 
\bigl[\rho:\pi_1(C\setminus c, x)\to U(E_x)\bigr]
$$
between
\begin{enumerate}
\item polystable vector bundles $E$ of rank $r$ and degree $d$
over a smooth projective curve $C$ 
with 2 marked points $x,c$, and
\item isomorphism classes of unitary representations
$\rho: \pi_1(C\setminus c, x)\to U(\c^r)$
of type $d/r$.
\end{enumerate}
\item $NS$
depends real analytically on $(C,c,x,E)$.
\item  The
fiber map $F\to F_x$ gives 
 a one--to--one correspondence between
\begin{enumerate}
\item direct summands of $E^{\otimes m}\otimes (E^*)^{\otimes n} $, and 
\item $\pi_1(C\setminus c, x)$ invariant subspaces of 
$E_x^{\otimes m}\otimes (E_x^*)^{\otimes n}$.
\end{enumerate}
\end{enumerate}
(This is stated in \cite{nar-ses} for $0\leq d <r$.
In the general case, we twist $E$ by a suitable
$\o_C(mc)$ and then apply \cite{nar-ses}.)

Because of the artificial role of the point $c$,
one has to be careful in taking determinants.
The representation $\det \rho$ corresponds to the degree 0
line bundle $\o_C(-d[c])\otimes \det E$.
\end{say}

\begin{defn} Let $E$ be a stable vector bundle
on a smooth projective curve $C$ and $x,c\in C$ closed points.
The Zariski closure of the image of the
Narasimhan--Seshadri representation
$\rho_c:\pi_1(C\setminus c,x)\to GL(E_x)$ is called
the {\it algebraic monodromy group} of $E$ at $(C,x,c)$
and it is denoted by $M_x(E,C,c)$.
Note that  $M_x(E,C,c)$ is reductive since it is the
Zariski closure of a subgroup of a unitary group.

$M_x(E,C,c)$ depends on the point $c$ but only slightly.
Choosing a different $c$ corresponds to tensoring
$E$ with a different line bundle, which
changes the representation by a character $\pi_1(C\setminus c,x)\to \c^*$.

As we see in (\ref{generic.monod.lem}),
for very general $c\in C$ we get the same
$M_x(E,C,c)$. We denote this common group
by $M_x(E,C)$. 
\end{defn}

\begin{lem} \label{scalars.claim} If $\det E$ is torsion in $\pic C$ then
the image of $\det: M_x(E,C)\to \c^*$ is torsion.
Otherwise $M_x(E,C)$ contains the scalars $\c^*\subset GL(E_x)$.
\end{lem}

Proof. As we noted above, 
$$
E_{\det \rho_c}\cong \o_C(-(\deg E)[c])\otimes \det E.
$$
If $\deg E\neq 0$, then
$E_{\det \rho_c}$ is a nonconstant family of
degree zero line bundles on $C$, hence its general member
is not torsion in $\pic C$. Thus in this case
$\det: M_x(E,C)\to \c^*$ is surjective.

If $\deg E=0$ then $E_{\det \rho_c}\cong \det E$
is constant. Thus $\det: M_x(E,C)\to \c^*$ is surjective
iff $\det E$ is not torsion in $\pic C$.

 Since $M_x(E,C)$ is reductive, 
we see that $\det: M_x(E,C)\to \c^*$ is surjective iff
the center of $M_x(E,C)$ is positive dimensional.

If $E$ is stable then $M_x(E,C,c)$ acts
irreducibly on $E_x$, 
and so the center consists of scalars only.
Thus we conclude that if $\det E$ is not torsion in $\pic (X)$
then
the scalars are contained in $M_x(E,C)$.
In general,
it is easy to see that 
$M_x(E,C)=\c^*\cdot M_x(E,C,c)$ for any $c\in C\setminus x$
if $\det E$ is not torsion in $\pic (X)$. \qed
\medskip

 The Narasimhan--Seshadri representations
$\rho$ vary  real analytically  with $(C,c,x,E)$
but the variation is definitely not complex analytic.
So it is not even clear that the  groups
$M_x(E,C)$ should vary algebraically in any sense.
Nonetheless, the situation turns out to be quite reasonable.

\begin{lem}\label{generic.monod.lem}
 Let $g:U\to V$ be a flat family of smooth projective curves
with  sections $s_x,s_c:V\to U$. Let $E\to U$ be a vector bundle
of rank $r$ 
such that $E|_{U_v}$ is polystable for every $v\in V$.
For every $v\in V$ let 
$$
\rho_v:\pi_1(U_v\setminus s_c(v),s_x(v))\to U(E_{s_x(v)})
$$
be the corresponding  Narasimhan--Seshadri  representation and let
$M_v\subset GL(E_{s_x(v)})$ be the Zariski closure of its image.

Then there is an open set $V^0\subset V$ and a
flat, reductive group scheme $G\subset GL(s_x^*E)\to V^0$ such that
$M_v=G_v$ for  very general $v\in V^0$.
(That is, for all $v$ in the complement of countably many
subvarieties of $V^0$.)
\end{lem}

\begin{rem}
By \cite[3.1]{rich},  the fibers of a flat, reductive group scheme
are conjugate to each other. The conjugacy class of the fibers
$G_v\subset GL(\c^r)$ is called the {\it generic monodromy}
group of $E$ on $U/V$. Note that while the
monodromy groups $M_x(E,C)$ are subgroups of $GL(E_x)$,
the generic monodromy
group is only a conjugacy class of subgroups.

In most cases $M_v=G_v$ for  every  $v\in V^0$, but there are
many exceptions. The simplest case is when $V=C$ is an elliptic curve,
$U=C\times C$ and $E$ is the universal degree 0 line bundle.

Then $M_c=\c^*$ if $c\in C$ is not torsion in $C$
but $M_c=\mu_n$, the groups of $n$th roots of unity,
if $c\in C$ is $n$-torsion.
\end{rem}

\begin{say}[Proof of (\ref{generic.monod.lem})]
 Let $W$ be a vector space of dimension $r$.
The general orbit of $GL(W)$ on
$\bigl(W^r+{\det}^{-1}W\bigr)^*$ is closed, hence the same holds for
any closed subgroup of $GL(W)$. We can thus recover the stable orbits of $G$,
and hence $G$ itself, as the general fibers of the rational map
$$
h_W:\bigl(W^r+{\det}^{-1}W\bigr)^* \map \spec
 \textstyle{\sum_{m\geq 0}} \bigl(S^{mr}(W^r)\otimes {\det}^{-m}W\bigr)^G.
$$
Correspondingly, if $E\to C$ is a rank $r$  vector bundle corresponding to
a unitary representation  
$\rho:\pi_1(C\setminus c,x)\to U(E_x)$, then we can recover the
Zariski closure of $\im\rho$  
from the general fibers of the rational map
$$
\begin{array}{ccl}
h_C:\bigl(E_x^r+{\det}^{-1}E_x\bigr)^* &\map 
&\spec
 \textstyle{\sum_{m\geq 0}} \bigl(S^{mr}(E_x^r)\otimes {\det}^{-m}E_x\bigr)^G\\
&& \qquad\qquad\qquad\downarrow\cong\\
&& \spec
 \textstyle{\sum_{m\geq 0}} H^0\bigl(C,S^{mr}(E^r)\otimes {\det}^{-m}E\bigr).
\end{array}
$$

Let us now apply this to our family $g:U\to V$.
Then we get a rational map
$$
h_V:\bigl(s_x^*E^r+{\det}^{-1}s_x^*E\bigr)^* \map
\textstyle{\sum_{m\geq 0}} g_*\bigl(S^{mr}(E^r)\otimes {\det}^{-m}E\bigr).
$$
Each of the sheaves
$$
g_*\bigl(S^{mr}(E^r)\otimes {\det}^{-m}E\bigr)
$$
commutes with base change over an open set $V_m\subset V$,
but these open sets may depend on $m$.
By the above remarks, for every point $v\in  \cap_{m\geq 1}V_m$,
a general fiber of $h_V$ above $v$ is 
the Zariski closure of the unitary representation
$\rho_v$.

Over the generic point $v_{gen}\in V$ we get a
reductive group scheme $G_{gen}\subset GL(E_{s(v_{gen})})$
which extends to a reductive group scheme $G\subset GL(s^*E|_{V^0})\to V^0$
over a suitable open set $V^0$.

The very general points in the lemma will be, by definition, the points in
the intersection $\cap_{m\geq 0}V_m$.
\qed
\end{say}

By taking the closure of $G$ in $GL(s^*E)$, we obtain an
open subset $V^*\subset V$ such that
\begin{enumerate}
\item the closure of $G$ in $GL(s^*E|_{V^*})$ is a flat
group scheme (but possibly not reductive), and
\item $V\setminus V^*$ has codimension $\geq 2$ in $V$.
\end{enumerate}

\begin{lem}\label{mon.fam.1d.lem} Notation as above.
For every $v\in V^*$, 
\begin{enumerate}
\item $M_v\subset  G^*_v$, 
\item $M_v$ is conjugate to a subgroup of the generic monodromy group, and
\item if  $\dim M_v = \dim G^*_v$ then in fact $M_v=G^*_v$.
\end{enumerate}
\end{lem}

Proof.  $U\to V$ is topologically a product in a Euclidean neighborhood
of $v\in W\subset V^*$, thus we can think of the family of representations
$\rho_v$ as a continuous map
$$
\rho: W\times \pi_1(U_v\setminus c(v),x(v))\to GL(\c^r).
$$
By (\ref{generic.monod.lem}), for very general $w\in W$,
$\rho(\{w\}\times \pi_1(U_v\setminus c(v),x(v))\subset G^*_w$,
hence, by continuity,
$\rho(\{v\}\times \pi_1(U_v\setminus c(v),x(v))\subset G^*_v$,
which proves (1).

Since $M_v$ is reductive, by \cite[3.1]{rich}, it is
conjugate to a subgroup of  $G^*_w$ 
for $w$ near $v$, hence to a subgroup of the generic monodromy group.

Finally, if $\dim M_v = \dim G^*_v$, then the connected component
of $G^*_v$ is the same as the connected component
of $M_v$, hence $G^*_v$ is reductive and again by
\cite[3.1]{rich}, it is
conjugate to a subgroup of the generic monodromy group.
Since $\rho(\{w\}\times \pi_1(U_v\setminus c(v),x(v))$ has points in every
connected component of $G^*_w$,
 by continuity the same holds for
$\rho(\{v\}\times \pi_1(U_v\setminus c(v),x(v))$.
Thus in fact $M_v=G^*_v$.\qed

\section{Holonomy groups}

Let $X$ be a normal, projective variety of dimension $d$
 with an  ample divisor
$H$.  A  curve $C\subset X$ is called a 
{\it complete intersection} (or  CI)  curve of type 
$(a_1,\dots,a_{d-1})$ if $C$ is a (scheme theoretic)
 intersection of $(d-1)$  divisors  $D_i\in |a_iH|$. 
We say that  $C\subset X$ is  a {\it general} 
CI  curve of type 
$(a_1,\dots,a_{d-1})$
if the  divisors  $D_i\in |a_iH|$ are all general.

If a smooth point $x\in X$ is fixed then
a general  CI curve of type $(a_1,\dots,a_{d-1})$ through $x$ is 
an intersection of $(d-1)$ general divisors  $D_i\in |a_iH|$,
each passing through $x$.

Let $E$ be a reflexive sheaf on $X$ such that $E$ is
$\mu$-stable with respect to $H$. By \cite{meh-ram1}
this is equivalent to assuming that
$E|_C$ is a stable vector bundle for a
general CI curve $C$ of type $(a_1,\dots,a_{d-1})$ for $a_i\gg 1$.

If $E$ is locally free at the points $x_1,\dots,x_s$, then
this is also equivalent to assuming that
$E|_C$ is a stable vector bundle for a
general CI curve $C$ of type $(a_1,\dots,a_{d-1})$ 
passing through the points $x_1,\dots,x_s$ for $a_i\gg 1$.
(While this stronger form of \cite{meh-ram1}
is not stated in the literature, it is easy to modify the
proofs to cover this more general case.)

\begin{defn}\label{hol.defn}
 Let $X$ be a normal, projective variety of dimension $n$
 with an  ample divisor
$H$ and $E$ a reflexive sheaf on $X$ such that $E$ is
$\mu$-stable with respect to $H$.
Assume that $E$ is locally free at $x$.

Let $B\subset X$ be the set of points where either
$X$ is singular or $E$ is not locally free.
Then $B$ has codimension at least 2 in $X$.
This implies that all general CI curves are contained in $X\setminus B$
and there is a one--to--one correspondence between
saturated subsheaves of the reflexive hull of 
$E^{\otimes m}\otimes (E^*)^{\otimes n}$
and saturated subsheaves of  
$E^{\otimes m}\otimes (E^*)^{\otimes n}|_{X\setminus B}$.

The {\it holonomy group} of $E$ at $x$
is the unique smallest subgroup  $H_x(E)\subset GL(E_x)$ such that:

For every smooth, pointed, projective curve $(D,d,y)$  
and every morphism $g:D\to X$ such that $g(y)=x$, $E$ is locally free
along $g(D)$ and
$g^*E$ is polystable, 
 the image of the Narasimhan--Seshadri 
 representation of $\pi_1(D\setminus d,y)$ is contained in 
$H_x(E)\subset GL(E_x)=GL((g^*E)_y)$.
\end{defn}

\begin{thm} \label{hol.thm}
Notation and assumptions as in (\ref{hol.defn}).
\begin{enumerate}
\item 
 Let $C$ be a very general CI curve of type $(a_1,\dots,a_{d-1})$
 through $x$ for $a_i\gg 1$.
Then the image of the  Narasimhan--Seshadri 
 representation of $\pi_1(C\setminus c,x)$ is Zariski dense in $H_x(E)$.
In particular, $H_x(E)$ is reductive.
\item For every $m,n$, the
fiber map $F\to F_x$  gives 
 a one--to--one correspondence between
direct summands of the reflexive hull of
$E^{\otimes m}\otimes (E^*)^{\otimes n} $ and 
$H_x(E)$ invariant subspaces of $E_x^{\otimes m}\otimes (E_x^*)^{\otimes n}$.
\item  The conjugacy class of $H_x(E)$ is the smallest  
reductive   conjugacy class   $G$
such that the structure group of $E$ can be reduced to $G$.
\end{enumerate}
\end{thm}

\begin{rem} 
For every curve $C$, the image of the
unitary representation of $\pi_1(C\setminus c,x)$ is contained in 
a  maximal compact subgroup of $H_x(E)$.
While $H_x(E)$ is well defined as a subgroup of $GL(E_x)$,
we do not claim that this  maximal compact subgroup of $H_x(E)$
is independent of $C$.
Most likely the opposite holds: 
 the maximal compact  subgroup  
 is independent of $C$ iff $E$ is a flat vector bundle
on $X\setminus \sing X$.
\end{rem}

\begin{say}[Proof of (\ref{hol.thm})]

Fix $(a_1,\dots,a_{d-1})$ such that $E|_C$ is stable
for a general CI curve $C$  of type
$(a_1,\dots,a_{d-1})$. By (\ref{generic.monod.lem}),
the conjugacy class of
$M_x(E,C,c)\subset GL(E_x)$  is independent of $C$ for very general
 $C$  of type
$(a_1,\dots,a_{d-1})$ and $c\in C$.
Denote this conjugacy class by
$M_x(a_1,\dots,a_{d-1})$.
First we show that these 
conjugacy classes
$M_x(a_1,\dots,a_{d-1})$ stabilize.

\begin{lem} 
There is a conjugacy class $M_x$ of subgroups of $GL(E_x)$ such that
if the $a_i$ are sufficiently divisible
then $M_x(a_1,\dots,a_{d-1})=M_x$.
\end{lem}

Proof. Fix a very general CI curve $C$  of type
$(a_1,\dots,a_{d-1})$ such that $E|_C$ is stable. We compare the
monodromy group $M_x(E,C,c)$ with the 
monodromy group $M_x(E,C_k,c_k)$ where $C_k$ is a very general
CI curve   of type
$ka_1,a_2,\dots,a_{d-1}$.

The divisors $D_2,\dots,D_{d-1}$ do not need changing,
so we may assume that $\dim X=2$.
Then $C$ is defined by a section $u\in H^0(X,\o_X(a_1H))$.
Choose a general  $v\in H^0(X,\o_X(ka_1H))$ vanishing at $x$
and consider the family of curves
$C'_t:=(u^k+t^kv=0)$. The general member is a
CI curve $C'_t$  of type
$(ka_1,a_2,\dots,a_{d-1})$ through $x$.

Note that $\supp C'_0=C$ 
but $C'_0$ has multiplicity $k$ along $C$.
The family is not normal along $C_0$ and we can normalize it by 
introducing the
new variable $u/t$. We then get a family of curves
$C_t$ such that $C_t=C'_t$ for $t\neq 0$ and
$C_0$ is a smooth curve, which is a degree $k$ cyclic cover  
$g:C_0\to C$ ramified 
at the intersection points $(u=v=0)$.

Since $C_0\to C$ is totally ramified at $x$,
we see that $g_*:\pi_1(C_0\setminus c_0,x)\to \pi_1(C\setminus c,x)$
is surjective where $c_0\in C_0$ is any preimage of $c$.
In particular,
$$
M_x(g^*E,C_0,c_0)=M_x(E,C,c).
$$

We can apply  (\ref{mon.fam.1d.lem}) to the family $\{C_t\}$
to  conclude
that 
$$
\dim M_x(a_1,\dots,a_{d-1})\leq \dim M_x(ka_1,\dots,a_{d-1}),
$$
and if equality holds then
$M_x(a_1,\dots,a_{d-1})$ and $M_x(ka_1,\dots,a_{d-1})$ are conjugate.

Thus if we choose $(a_1,\dots,a_{d-1})$ such that
$\dim M_x(a_1,\dots,a_{d-1})$ is maximal, then
$M_x(a_1,\dots,a_{d-1})$ and $M_x(b_1,\dots,b_{d-1})$ are conjugate
whenever $a_i|b_i$ for every $i$.
\qed
\medskip

Choose $(a_1, \dots,a_{d-1})$ 
and a very general CI curve of type $(a_1, \dots,a_{d-1})$
through $x$ 
such that
\begin{enumerate}
\item $M_x(a_1, \dots,a_{d-1})=M_x$, and
\item every stable summand of $T(E)$ restricts to a stable
bundle on $C$.
\end{enumerate}

\begin{claim} With the above notation,  $H_x(E)=M_x$.
\end{claim}

Proof. 
$M_x\subset H_x(E)$ by definition.

By assumption $M_x$ is the stabilizer of a nonzero vector
 $w_x\in T(E)_x=E_x^{\otimes m}\otimes (E_x^*)^{\otimes n}$, 
thus it corresponds to a direct summand
$\o_C\cong W_C\into T(E)|_C$ which in turn gives a
direct summand
$\o_X\cong W_X\into T(E)$ by the second assumption.

Pick any 
 smooth pointed curve $(D,d,y)$  
and a map $g:D\to X\setminus B$ such that $g(y)=x$ and
$g^*E$ is polystable.
Then $\o_D\cong g^*W\subset g^*(T(E))$ is a
direct summand, hence the  
 Narasimhan--Seshadri 
 representation of $\pi_1(D\setminus d,y)$ 
in $g^*(T(E_x))=T(g^*E_y)$ fixes $w$.
The stabilizer of $w$ is exactly $M_x$, hence
the monodromy group of $g^*E$ is contained in $M_x$.
Since this holds for any $(D,d,y)$,
we see that $H_x(E)=M_x$. \qed

\begin{claim} The stabilizer of 
$W|_{X\setminus B}\into T(E|_{X\setminus B})$ in $GL(E|_{X\setminus B})$
is a reductive subgroup scheme  $H\subset GL(E|_{X\setminus B})$
 whose fibers are
in the conjugacy class $M_x$.
The structure group of $E|_{X\setminus B}$ can be reduced to 
a conjugacy class $G$
iff some group in $G$ contains $M_x$.
\end{claim}

Proof. By construction $H_x=M_x$ is reductive, hence
there is a largest open set $X^0\subset X$ such that
the fibers $H_v$ are
in the conjugacy class $M_x$ for every $v\in X^0$.
Thus the  structure group of $E|_{X^0}$ can  be reduced to $M_x$.

Pick a very general CI curve $C$ of type
$(a_1,\dots,a_{d-1})$ such that $E|_C$ is stable and 
$M_x(E,C)=M_x$.  The stabilizer of every point
of $W|_C$ is conjugate to $M_x(E,C)$, which shows that
$X^0$ contains $C$. This implies that
$X\setminus X^0$ has codimension $\geq 2$ in $X$.

By Hartogs' theorem, a rational map from a normal variety
to an affine variety which is defined outside a codimension two 
set is everywhere defined, thus
 the  structure group of $E|_{X\setminus B}$ also  reduces to $M_x$.

Conversely, if the  structure group of $E$ can be reduced to 
 the conjugacy class $G\subset GL(E)$,
then the  structure group of $E|_C$ can also be reduced to $G$,
hence some group in  the conjugacy class  $G$ contains $M_x$.\qed

\medskip

It remains to show that (\ref{hol.thm}.1) holds for
 $(a_1,\dots,a_{d-1})$ sufficiently large. (So far we
have established (\ref{hol.thm}.1)
only for $(a_1,\dots,a_{d-1})$ sufficiently divisible.)

Fix now $(a_1,\dots,a_{d-1})$ such that
$M_x(a_1,\dots,a_{d-1})=M_x$. We claim that in fact
$M_x(b_1,\dots,b_{d-1})=M_x$ for every 
$b_i\geq 2a_i$.

Indeed, assume the contrary.
By (\ref{mon.fam.1d.lem}) 
we know that $M_x(b_1,\dots,b_{d-1})$ is conjugate to
a subgroup of $M_x$. Thus if
they are not equal, then there are $m,n$ 
and a vector $v\in E_x^{\otimes m}\otimes (E_x^*)^{\otimes n}$
which is stabilized by $M_x(b_1,\dots,b_{d-1})$
but not by $M_x$.

Correspondingly, if $D$ is a
 very general CI curve   of type
$b_1,\dots,b_{d-1}$, then $v$ corresponds to a
direct summand
$W_D\subset E^{\otimes m}\otimes (E^*)^{\otimes n}|_D$ 
which can not be obtained as a restriction of a
 direct summand of
$E^{\otimes m}\otimes (E^*)^{\otimes n}$. 
Thus there is a stable direct summand
$F\subset E^{\otimes m}$ such that $F|_D$ is not stable.
By the already proved case of (\ref{hol.thm}.1),
we know that $F|_C$ is stable. 
Hence by \cite[Lem.7.2.10]{huy-leh}, 
$F|_D$ is also stable, a contradiction.
\qed
\end{say}

\section{Parabolic bundles}

We briefly recall the correspondence between the category of parabolic
bundles on $X$ and the category of $G$--bundles on a suitable cover.

Let $D$ be an effective divisor on $X$. For a coherent sheaf $E$ on
$X$ the image of $E\otimes_{{\mathcal O}_X} {\mathcal O}_X(-D)$ in
$E$ is  denoted by $E(-D)$. The following definition of parabolic
sheaves was introduced in \cite{MY}.

\begin{defn}\label{newdefn} Let $E$ be a torsion-free
  ${\mathcal O}_X$--coherent sheaf on $X$. A {\it quasi--parabolic}
  structure on $E$ over $D$ is a filtration by ${\mathcal
    O}_X$--coherent subsheaves
$$
E\, =\, F_1(E)\, \supset\, F_2(E)\, \supset\, \cdots
\,\supset\, F_l(E)\,\supset\, F_{l+1}(E)\,=\, E(-D)
$$
The integer $l$ is called the {\it length of the filtration}.
A {\it parabolic structure} is a quasi--parabolic structure,
as above, together with a system of {\it weights}
\[
0\, \leq\,
{\alpha}_1\, < \, {\alpha}_2 < \, \cdots \, < \, {\alpha}_{l-1} \, < \,
{\alpha}_l \, < \, 1
\]
where the weight ${\alpha}_i$ corresponds to the subsheaf $F_i(E)$.
\end{defn}

We shall denote the parabolic sheaf defined above datum by the triple
$(E,F_*,{\alpha}_*)$.  When there is no confusion it will be denoted
by $E_*$.

For a parabolic sheaf $(E,F_*, {\alpha}_*)$ define
the following filtration $\{E_t\}_{t\in \mathbb R}$ of coherent
sheaves on $X$ parameterized by $\mathbb R$:
$$
E_t \hspace{.1in} := \hspace{.1in} F_i(E)(-[t]D)
\eqno{(\ref{newdefn}.1)}
$$
where $[t]$ is the integral part of $t$
and ${\alpha}_{i-1} < t - [t]
\leq {\alpha}_i$, with the convention that ${\alpha}_0 = {\alpha}_l -1 $
and ${\alpha}_{l+1} = 1$.

A {\it homomorphism} from the parabolic sheaf
$(E, F_*, {\alpha}_*)$ to another parabolic sheaf
$(E', F'_*, {\alpha}'_*)$ is a homomorphism from $E$ to $E'$
which sends any subsheaf $E_t$ into $E'_t$, where $t \in [0,1]$
and the filtration are as above.

If the underlying sheaf $E$ is locally free then $E_*$ will be called
a parabolic vector bundle. {\it In this section, all parabolic sheaves
  will be assumed to be parabolic vector bundles.}

We have the following equivalent definition:
\begin{defn}\label{olddefn} Let $X$ be a normal, projective variety and $D$
  an effective divisor.  A {\it quasi--parabolic filtration} on a
  sheaf $E$ is a filtration by subsheaves of the restriction $E|_{D}$
  of the sheaf $E$ to the parabolic divisor $D$:
\[
E|_{D} = {\mathcal F}^1_{D}(E) \supset {\mathcal F}^2_{D}(E) \supset
\ldots \supset {\mathcal F}^l_{D}(E) \supset {\mathcal F}^{l+1}_{D}(E) = 0
\] 
together with a system of weights 
\[
0\, \leq\,
{\alpha}_1\, < \, {\alpha}_2 < \, \cdots \, < \, {\alpha}_{l-1} \, < \,
{\alpha}_l \, < \, 1
\]
\end{defn}

We assume that the following conditions are satisfied:

\begin{enumerate}

\item  $X$ is smooth and $D$ is a
 divisors with normal crossings. In other words,
any parabolic divisor is assumed to be reduced, its  irreducible
components  are smooth and  the irreducible components
intersect transversally.  

\item All the parabolic weights are 
    rational numbers. 
  
\item On each component of
  the parabolic divisor the filtration is given by  subbundles.
\end{enumerate}

Consider the decomposition 

\[
D = \sum_{i=1}^{n} D_i
\]
Let $E$ be a vector bundle on $X$. For $1 \leq i \leq n$, let
\[
E|_{D_i} = {\mathcal F_i}^1 \supset {\mathcal F_i}^2 \supset \ldots
\supset {\mathcal F_i}^{l_i} \supset {\mathcal F_i}^{l_i + 1} = 0
\] 
with $l_i \geq 1$, be a filtration of subbundles on $D_i$. Suppose
that we are given a string of numbers $\alpha^i_j$, with $1 \leq j
\leq l_i + 1$, satisfying the following:
\[
0\, \leq\, \alpha^i_1\, < \, \alpha^i_2 < \, \cdots \, < \,
\alpha^i_{l_i} \, < \, \alpha^i_{l_i + 1} = 1
\]
Then we can construct a parabolic structure on $E$ as follows: Define
the coherent subsheaves ${\overline{F_i}}^j$ of $E$, where $1
\leq j \leq l_i$ by the following short exact sequence:
\[
0 \to {\overline{F_i}}^j \to E \to
(E|_{D_i})/{\mathcal F_i}^j \to 0
\]
For $1 \leq i \leq n$ and $0 \leq t < 1$, let
\[
l(t,i) := \min\bigl[j |~j\in \{1, \ldots , l_i + 1 \} ~\&~~ \alpha^i_j
\geq t \bigr]
\]
Define 
\[
E_t = \cap_{i = 1}^{n} {\overline{F_i}}^{l(t,i)}
\subset E
\]
The filtration $\{E_t\}$ defines a parabolic structure on $E$ and any
parabolic structure on $E$ with $D$ as parabolic divisor arises this
way.

We  denote the entire parabolic datum by $(E, F_*, \alpha_*)$ or
simply by $E_*$ when the context is clear. If the underlying sheaf $E$
is locally free then $E_*$ is called a parabolic vector bundle.

Let ${\rm PVect}(X,D)$ denote the category whose objects are parabolic
vector bundles over $X$ with parabolic structure over the divisor $D$
satisfying the above three conditions, and the morphisms of the
category are homomorphisms of parabolic vector bundles (cf.\ for
example \cite{biswas}).

The direct sum of two vector bundles with parabolic structures has an
obvious parabolic structure and ${\rm PVect}(X,D)$ is closed under the
operation of taking direct sum. We remark that the category ${\rm
  PVect}(X,D)$ is an additive tensor category with the direct sum and
the parabolic tensor product operation. It is straight--forward to
check that ${\rm PVect}(X,D)$ is also closed under the operation of
taking the parabolic dual defined in \cite{biswas} or \cite{yokogawa}.

For an integer $N\geq 2$, let ${\rm PVect}(X,D,N) \, \subseteq \, {\rm
  PVect}(X,D)$ denote the subcategory consisting of all parabolic
vector bundles all of whose parabolic weights are multiples of $1/N$.
It is straight--forward to check that ${\rm PVect}(X,D,N)$ is closed
under all the above operations, namely parabolic tensor product,
direct sum and taking the parabolic dual.

\begin{say}[The covering construction]
Let $X$ be a smooth projective variety and $D$
an effective simple normal crossing divisor.
The Covering Lemma of Kawamata
 \cite[Thm.17]{K} says
that there is a connected smooth projective
variety $Y$  and a Galois covering morphism
$$
p \,: \hspace{.1in} Y\hspace{.1in} \longrightarrow
\hspace{.1in} X 
$$
such that the reduced divisor $\tD:= \,({p}^*D)_{red}$
is a normal crossing divisor on $Y$ and furthermore,
${p}^*D_i= k_iN\cdot({p}^*D_{i})_{red}$,
where the $k_i$
are positive integers. Let $G$ denote the Galois group
for the covering map $p$.
\end{say}

\begin{defn}[The category of $G$--bundles]\label{jumbo}
Let $G \, \subseteq \, \mbox{Aut}(Y)$ be a finite subgroup of the
group of automorphisms of a connected smooth projective variety
$Y$. The natural action of $G$ on $Y$ is encoded in a morphism
$$
{\mu} \,:\hspace{.1in} G \times Y \hspace{.1in}
\longrightarrow \hspace{.1in} Y
$$

Let ${\rm Vect}_{G}(Y)$ denote the category of all $G$--linearized
vector bundles on $Y$.  The isotropy group of any point $y \in Y$, for
the action of $G$ on $Y$, will be denoted by ${G}_y$.

Let ${\rm Vect}^D_{G}(Y,N)$ denote the subcategory of ${\rm
  Vect}_{G}(Y)$ consisting of all $G$--linearized vector bundles
$W$ over $Y$ satisfying the following three conditions:

\begin{enumerate}
  
\item for a general point $y$ of an irreducible component of
  $(p^*D_i)_{red}$, the isotropy subgroup ${G}_y$ is 
    cyclic of order $|{G}_y| = n_y$ which is a divisor of $N$; the
  action of the isotropy group ${G}_y$ on the fiber $W_y$ is of order
  $N$, which is equivalent to the condition that for any $g \in
  {G}_y$, the action of $g^{N}$ on $W_y$ is the trivial action;
  
\item The action is given by a representation $\rho_y$ of
  ${G}_y$ given by a block diagonal matrix
$$
\rho_y(\zeta)=\diag\bigl({z}^{{\alpha}_1}I_1,\dots,{z}^{{\alpha}_l}I_l\bigr)
$$
where
 $\zeta$ is a generator of the group ${G}_y$ and whose order
$n_y$ divides $N$, 
 $\alpha_i = \frac{m_j}{N}$,
 $I_j$ is the identity matrix of order $r_j$, where $r_j$ is the
  multiplicity of the weight $\alpha_j$, $z$ is an $n_y$-th root of unity
and 
 $0 \leq m_1 < m_2 < ...< m_l \leq N - 1 $.

\item For a general point $y$ of an irreducible component of a
  ramification divisor for $p$ not contained in $(p^*D)_{red}$,
  the action of ${G}_y$ on $W_y$ is the trivial action.

\end{enumerate}

Following Seshadri \cite[p.161]{seshadri} we call the $G$--bundles
$E$ in ${\rm Vect}^D_{G}(Y,N)$ {\it bundles of fixed local orbifold
  type $\tau$}. 

We remark that this definition of $G$--bundles of fixed local
type easily extends to $G$--torsion--free sheaves since the local
action is specified only at the generic points of the parabolic
divisor. 

We note that ${\rm Vect}^D_{G}(Y,N)$ is also an additive tensor
category.
\end{defn}

\begin{say}[Parabolic bundles and $G$--bundles]\label{parabgamma}

In \cite{biswas} an identification between the objects of ${\rm
  PVect}(X,D,N)$ and the objects of ${\rm Vect}^D_{G}(Y,N)$ has been
constructed. Given a $G$--homomorphism between two $G$--linearized
vector bundles, there is a naturally associated homomorphisms between
the corresponding vector bundles, and this identifies, in a bijective
fashion, the space of all $G$--homomorphisms between two objects of
${\rm Vect}^D_{G}(Y,N)$ and the space of all homomorphisms between
the corresponding objects of ${\rm PVect}(X,D,N)$.  An equivalence
between the two additive tensor categories, namely ${\rm
  PVect}(X,D,N)$ and ${\rm Vect}^D_{G}(Y,N)$, is obtained this way.

We observe that an earlier assertion that the parabolic tensor product
operation enjoys all the abstract properties of the usual tensor
product operation of vector bundles, is a consequence of the fact that
the above equivalence of categories indeed preserves the tensor
product operation.

The above equivalence of categories has the further property that it
takes the parabolic dual of a parabolic vector bundle to the usual
dual of the corresponding $G$--linearized vector bundle.
\end{say}

\begin{defn}[Stable parabolic bundles]
The definition of parabolic semistable and parabolic polystable vector
bundles is given in Maruyama-Yokogawa \cite{MY} and Mehta-Seshadri
\cite{MS}. Given an ample divisor $H$, the
{\it parabolic degree} of a parabolic bundle $E_*$ is defined by

\begin{equation}\label{parabdeg} 
\pardeg(E_*) := \int^{1}_{0} \deg(E_t) dt + r\cdot \deg(D) 
\end{equation}
where $E_t$ comes from the filtration defined in
(\ref{newdefn}.1). There is a natural notion of parabolic
subsheaf and given any subsheaf of $E$ there is a canonical parabolic
structure that can be given to this subsheaf.  (cf \cite{MY}
\cite{biswas} for details)

A parabolic bundle is called {\it stable} (resp. {\it semistable}) if
for any proper nonzero coherent parabolic subsheaf $V_*$ of $E_*$ with
$0 < rank(V_*) < rank(E_*)$,, with $E/V$ being torsion free, the
following inequality is valid:
$$
{{\pardeg V}\over {\rank V}} \, < \,
{{\pardeg E}\over {\rank E}} \hspace{.3in}
\qtq{resp.} {{\pardeg V}\over {\rank V}} \, \leq \,
{{\pardeg E}\over {\rank E}}.
$$
\end{defn}

\begin{rem}\label{oldpardeg} If we work with the definition given in 
  \eqref{olddefn}, then we have the following expression for parabolic
  degree of $E_*$ which is along the lines of \cite{MS}. Define:
\[
\weight(E_*) := \sum_{i,j} \alpha^i_j \bigl(c_1({\mathcal F_i}^j(E)) \cdot
H^{n-1} - c_1({\mathcal F_i}^{j+1}(E)) \cdot H^{n-1} \bigr)
\]
Using the fact that on the divisor $D$, $c_1(F) = rk(F) D$, we have
the following expression:
\[
\weight(E_*) := \sum_{i,j} \alpha^i_j \bigl[ \rank({\mathcal F_i}^j(E)) -
\rank({\mathcal F_i}^{j+1}(E))\bigr] \bigl(D_i\cdot H^{n-1} \bigr)
\]
Then it is not hard to check that the parabolic degree of $E_*$ is
given by:
\[
\pardeg(E_*) = \deg(E) + \weight(E_*)
\]
\end{rem}

\begin{defn}[Stable $G$-bundles]
  A $G$-linearized vector bundle $V'$ over $Y$ is called $(G,
  \mu)$-{\it stable} (resp.\ $(G, \mu)$-{\it semistable}) if for any
  proper nonzero coherent subsheaf $F' \subset V'$, invariant under
  the action of $G$ and with $V'/F'$ being torsion free, the following
  inequality is valid:
$$
{{\deg F'}\over {\rank F'}} \, < \,
{{\deg V'}\over {\rank V'}} \hspace{.3in}
\qtq{resp.} {{\deg F'}\over {\rank F'}} \, \leq \,
{{\deg V'}\over {\rank V'}}.
$$
The $G$-linearized vector bundle $V'$ is called
$G$-{\it polystable} if it is a direct sum of
$G$-stable vector bundles of same slope, where, as usual,
$\mbox{slope}:=\mbox{degree}/
\mbox{rank}$. 
\end{defn}

\begin{rem}\label{biswas}
$G$-invariant subsheaves of $V'$ are in one-to-one
correspondence with the subsheaves of the parabolic vector bundle
corresponding to $V'$, and furthermore, the degree of a
$G$-invariant subsheaf is simply the order of $G$-times the parabolic degree
of the corresponding subsheaf with the induced parabolic structure
\cite{biswas}. 

It is immediate that  $V^*$ is $G$-semistable
 if and only if $V$ is so.

The above equivalence of categories between ${\rm
  PVect}(X,D,N)$ and ${\rm Vect}^D_{G}(Y,N)$ in fact identifies the
subcategory of {\it parabolic stable} bundles with the $G$--stable
bundles. This result, due to Biswas, generalizes the result of Seshadri
for parabolic bundles over curves (cf \cite{biswas}, \cite{seshadri}).
\end{rem}

\begin{prop} Let $E$ be a stable vector bundle on $X$ with $rank(E) = n$ and
  $deg(E) = q$ and such that $-n < q \leq 0$.  Then, for any smooth
  divisor $D \subset X$ such that $D \in |H|$, one can endow $E$ with
  a parabolic structure along $D$ such that $\pardeg(E) = 0$ and $E$
  is parabolic stable with this structure.
\end{prop}

Proof. Let $p : Y \to X$ be a Kawamata cover of $X$ with Galois
group $G$ and ramification index along $D$ to be the integer $n$.
Define $\tilde D := (p^{*}(D))_{red}$ so that $p^{*}(D) = n \cdot
(p^{*}(D))_{red}$. Further, in the notation of (\ref{jumbo}), the
weight $\alpha$ attached to the action of the isotropy $G_{y}$ at a
general point $y \in {\tilde D}$ is given by $\alpha = \frac{-q}{n}$.

Since $\tilde D$ is invariant under the action of $G$, for any $k
\in {\bf Z}$, the line bundle ${\mathcal O}_Y(k {\tilde D})$ gets a
structure of a $G$--bundle.

Define ${\bf L} = {\mathcal O}_Y(-q \cdot {\tilde D})$. Then ${\bf L}$
also gets a $G$--bundle structure. Now consider the $G$--bundle
$p^{*}(E)$ and let $W$ be the $G$--bundle (of type $\tau$ in the
notation of (\ref{jumbo}) defined by:
\[
W = p^{*}(E) \otimes_{{\mathcal O}_Y} {\bf L}
\]
It is easy to see that $p^{G}_*(W) = E$. Further, $E$ realised as
the invariant direct image of $W$ gets a natural parabolic structure,
called the {\it special} parabolic structure where the flag has only
two terms 
\[
E|_{D} = {\mathcal F}^1_{D}(E) \supset {\mathcal F}^2_{D}(E) = 0
\] 
with weight $\alpha = \frac{-q}{n}$.

The parabolic degree of $E$ with this structure
is given by:
\[
\pardeg(E) = \deg(E) + n \cdot \alpha = \deg(E) - q = 0.
\]
We observe that for any subbundle $V \subset E$ with $\rank(V) = r$,
there is a unique way of defining the induced {\it special} parabolic
structure on $V$ and $\pardeg(V) = \deg(V) + r \cdot \alpha = \deg(V) + r
\cdot \frac{-q}{n}$.  Hence,
\[
\frac{\pardeg(V)}{r} = \frac{\deg(V)}{r} + \frac{-q}{n} <
\frac{\pardeg(E)}{n} = 0
\]
since $E$ is stable. Thus, we conclude that $E$ is
{\it parabolic stable} with this parabolic structure. We also note
that by the correspondence between parabolic stable bundles on $X$ and
$(G, \mu)$--stable bundles on $Y$ (Remark \ref{biswas}), the
$G$--bundle $W$ is $(G, \mu)$--stable. \qed
\medskip

\begin{rem} This proposition can be seen in the more general context of
  parabolic bundles. Let $E_* \in Vect(X,D)$ with $\pardeg(E_*) \neq
  0$. Then there exists a parabolic bundle $E'_* \in Vect(X,D')$,
  where $D \subset D'$ and $D'$ has more components ${\Delta_j}$
  meeting $D_i$ with simple normal crossing singularities, such that
$$
\pardeg(E'_*) = 0
$$
and further, the  forgetful functor  
\[
Vect(X,D') \to Vect(X,D), ~~~E'_* \to E_*
\]
is fully faithful and preserves parabolic semistability and parabolic
stability. 

To see this, we define $E'_* = (E',F'_*,\alpha'_*)$ as follows: 

Assume that $\pardeg(E_*) < 0$. This is always possible to achieve by
twisting with a line bundle. Let integers $m_j > 0$ be so chosen, such
that for rational numbers $0 \leq \beta_j < 1$, we have the equality
\[
(\sum (m_j \beta_j (\Delta_j \cdot H^{n-1})) =
\frac{-\pardeg(E_*)}{\rank(E)}
\]

Let $D' = \sum_{i = 1}^{n} D_i + \sum_{j=1}^{m} {\Delta_j}$ and $E' =
E$. Define the filtration as follows:
\[
E'|_{\Delta_j} = {\mathcal F_{0}}^1(E') \supset  {\mathcal F_{0}}^2(E') = (0)  
\]
with a single weight $\beta_j$ for each $1 \leq j \leq m$. That is,
$(\alpha'_*) = (\alpha_*) \cup (\beta_*)$.  Clearly, $\weight(E'_*) =
\weight(E_*) + \rank(E) \bigl[\sum (m_j \beta_j (\Delta_j \cdot H^{n-1})
 \bigr]$.
Hence, $\pardeg(E'_*) = \pardeg(E_*) + \rank(E) \bigl[\sum (m_j \beta_j
(\Delta_j \cdot H^{n-1})\bigr] = 0$.

\end{rem}
 
The above Proposition is used to extend our theory of holonomy to
parabolic stable bundles.  Before we do that, we need to prove the
Mehta--Ramanathan restriction theorem for $G$-torsion free sheaves. By
the equivalence of categories between $G$--bundles and parabolic
bundles, we get a Mehta--Ramanathan--type restriction theorem for
parabolic bundles.

\begin{thm}[The $G$--Mehta-Ramanathan theorem] Let $E$ be a
$(G,\mu)$--semistable (resp.\ stable) $G$--torsion free sheaf on a
normal projective $G$--variety. Then the restriction $E|_{C}$  to a general
complete intersection $G$--curve $C$ of large degree (with respect
to the pull-back line bundle $p^{*}(H)$ ) is $(G,\mu)$--semistable
(resp.\ stable).  
\end{thm}

Proof.  Since $(G,\mu)$--semistability for $G$--sheaves is equivalent
to the semistability of the underlying sheaf, the non-trivial case is
that of stability. The proof follows from the following observations:

\begin{enumerate}
\item $E$ is $(G,\mu)$--stable iff $E$ is polystable
and $\Hom(E,E)^G$ is 1--dimensional.

Indeed, we noted that $E$ is semistable. If $E$ is not 
polystable then it has a nontrivial 
socle $F\subset E$ with $\mu(F) = \mu(E)$ which is invariant under all
  the automorphisms of $E$, in particular invariant under the group
  $G$ (cf.\ \cite[1.5.5]{huy-leh}). This contradicts the $G$ stability of $E$.

\item  By the orbifold version of the Enriques-Severi 
theorem, for
  sufficiently high degree complete intersection $G$--curve $C$,
 $\Hom_X(E,E)=\Hom_C(E|_C,E|_C)$ 
and so  $\Hom_X(E,E)^G=\Hom_C(E|_C,E|_C)^G$. 

\item Finally,  by the  restriction theorem of
  Bogomolov (cf.\ \cite[Sec.7.3]{huy-leh}),
 for {\it every} complete intersection curve $C$ in
  the linear system $|m H|$ (the number $m$ being effectively
  determined), the restriction $E|_{C}$ is polystable.
Thus we can work with general high degree complete intersection $G$--curves
rather than with general complete intersection curves. \qed
\end{enumerate}

We can now define the holonomy groups of  parabolic stable bundles.

\begin{defn}[Holonomy groups of  parabolic stable bundles]
\label{parab.holom.sefd}
Let $E$ be a $(G,\mu)$--stable bundle on $Y$ of degree $0$.
  This corresponds to a parabolic
  bundle on $X$ of  parabolic degree $0$.
 Let $C_k$ be a general CI curve in
  $Y$ which is $G$--invariant.
  
The quotient  $C_k/{G} =: T_k$ is also a
   smooth projective curve in $X$. By choosing $C_k$ sufficiently
  general, one can make sure that the action of $G$ on $C_k$ is
  faithful and  we can realize the group $G$ as a
  quotient $\Gamma/{\Gamma}_o$, where ${\Gamma}_o = {\pi}_1(C_k)$ and
 $\Gamma$ acts
  properly discontinuously on the simply connected cover ${\tilde
    C_k}$ and  $T_k = {\tilde C_k}/{\Gamma}$. (The $\Gamma$--action on
  ${\tilde C_k}$ is not assumed to be free.)
  
 By the restriction theorem above, $E|_{C_k}$ is a
  $(G,\mu)$--stable bundle on $C_k$, hence  it comes from an irreducible
  unitary representation of the group $\rho: \Gamma \to GL(E_y)$, for $y
  \in C_k \subset Y$ a point away from the ramification locus. We note
  that an irreducible unitary representation of $\Gamma$ descends to a
  bundle on ${\tilde C_k}/{\Gamma}_o = C_k$ which comes with a
  $G$--action.

   Now by considering the map $p: C_k \to T_k$ and taking the
  invariant direct image $p_{*}^{G}(E|_{C_k})$ we get a bundle $F$
  which is  parabolic stable on $T_k$, with parabolic structure on
  $T_k \cap D$. Hence,  as above, the group $\Gamma$
  which acts on the simply connected curve ${\tilde C_k}$ with
  parabolic fixed points such that $p_{*}^{G}(E|_{C_k}) = F$ arises
  from a unitary representation $\rho : \Gamma \to GL(F_x)$, with $p(y) =
  x$.

    Let $y \in C_k \subset Y$ be a point away from the ramification
  locus. The arguments in Section 2
 now imply that the Zariski closure of $\im (\rho) = H_y$
  is well-defined and is the smallest reductive subgroup of $GL(E_y)$
  such that the $G$--bundle $E$ has a reduction of structure group
  to $H_y$  and the reduction is moreover $G$--invariant.
  Moreover, $H_y$ can be identified with the Zariski closure of the
  image of $\rho$ in $GL(F_x)$ where $p(y) = x$.
  
   By the categorical equivalence between $G$--bundles on $Y$ and
  parabolic bundles on $X$, it follows that the group $H_y = H_x$ is
 realized as the  holonomy group of the parabolic
    bundle $p^{G}_{*}(E)$ on $X$. This defines the holonomy group
  for all parabolic
    stable bundles  in the category ${\rm PVect}(X,D)$.
\end{defn}
  
\begin{rem} If the bundle arises as an irreducible representation of
$\pi_1(X \setminus D)$ then the resulting parabolic bundle will have all
parabolic Chern classes zero and this fits into the theme addressed by
Deligne in \cite{deligne}.\end{rem}

\section{Computing the holonomy group}

Given a stable vector bundle $E$, the computation of its
holonomy group seems quite hard in general.
The definition (\ref{hol.defn}) is practically
impossible to use. The method of Tannaka duality \cite{tannaka}
shows that one can determine the holonomy
once we know the decomposition of
$E^{\otimes m}\otimes (E^*)^{\otimes n}$ into direct summands for every $m,n$.
The observation of Larsen (which seems to be unpublished) is that
one can frequently characterize
a subgroup $G\subset GL(V)$ by knowing
the decomposition of the $G$-module
$V^{\otimes m}\otimes (V^*)^{\otimes n}$ for only a very few values of $m,n$.

Our aim is to translate this into geometric form
and give several examples illustrating the
principle (\ref{hsm.princ}).
Let us start with the general form of
(\ref{hol.simpconn.prop}).

\begin{lem} \label{hol.simpconn.lem}
Notation and assumptions as in (\ref{hol.defn}).
There is an \'etale cover
$\pi:U\to (X\setminus \sing X)$
with Galois group $H_x(E)/H_x(E)^0$
such that the holonomy group of
$\pi^*E$ is   $H_x(E)^0$, hence  connected.
\end{lem}

Proof. Let $B\subset X$ be the set of points where either
$X$ is singular or $E$ is not locally free.
To the vector bundle $E|_{X\setminus B}$
we can associate a principal $H_x(E)$-bundle $P\to (X\setminus B)$.
Then $U:=P/H_x(E)^0\to (X\setminus B)$ is an \'etale cover
with Galois group $H_x(E)/H_x(E)^0$.
Since $\pi^*P/H_x(E)^0\to U$ has a section, 
 the structure group of $\pi^*P$ can be further reduced to 
$H_x(E)^0$.\qed

The following result relates the holonomy groups to symmetric powers.

\begin{prop}\label{S^m.cond.pf}
Notation  as in (\ref{hol.defn}).
Assume that  $X\setminus B$ is  simply connected.
Then the following are equivalent:
\begin{enumerate}
\item The reflexive hull of $S^mE$ is indecomposable for some $m\geq 2$.
\item The reflexive hull of $S^mE$  is indecomposable for every $m\geq 2$.
\item The holonomy is one of the following: 
\begin{enumerate}
\item $SL(E_x)$ or $GL(E_x)$, 
\item $Sp(E_x)$ or $GSp(E_x)$ 
for a suitable nondegenerate  symplectic form 
on $E_x$ (and  $\rank E$ is even).
\end{enumerate}
\end{enumerate}
\end{prop}

Proof. Let ${\frak h}\subset {\frak gl}(E_x)$
denote the Lie algebra of $H_x(E)$.
A representation of $H_x(E)$ is indecomposable iff
the corresponding representation of ${\frak h}$ is indecomposable.
Thus (\ref{S^m.cond.pf}) is equivalent to the
corresponding statement about Lie algebras.
The latter is a special case of \cite[\S 17,Thm.1]{seitz}\qed
\medskip

The following is a key example in relating the
holonomy groups to geometric structures.

\begin{exmp}\label{s03.exmp} Let $E$ be a rank 3 bundle with 
holonomy group  $SO_3$ over a smooth projective variety $X$

We can also think of the standard representation of
$SO_3\cong PSL_2$ as the symmetric square
of the standard representation of $SL_2$.
Does this mean that every rank 3 bundle with
$SO_3$-holonomy can be written as the symmetric square
of a rank 2 bundle with $SL_2$-holonomy?

Principal $PSL_2$-bundles are classified
by $H^1_{et}(X,PSL_2)$.
The obstruction to lift to
a principal $SL_2$-bundle is
in $H^2_{et}(X,\mu_2)$, which is never zero.
(For a basic reference, see \cite[Sec.IV.2]{milne}.)

To put it in more concrete terms,
 observe that for any rank 2 bundle $F$,
the rank 3 bundle $S^2F\otimes {\det}^{-1}F$
has trivial determinant and $SO_3$-holonomy. 
If $\det F\cong L^2$ is the square of a line bundle,
then
$$
S^2F\otimes {\det}^{-1}F\cong S^2(F\otimes L^{-1}),
$$
but if $\det F$ is not the square of a line bundle,
then there does not seem to be any natural way to
write $S^2F\otimes {\det}^{-1}F$ as a symmetric square.

This is the obstruction in $\pic(X)/\pic(X)^2\subset H^2_{et}(X,\mu_2)$
that we detected earlier.

This suggests that it is 
easier to lift a $PSL_2$-bundle to a $GL_2$-bundle
tensored with a line bundle than to an $SL_2$-bundle.
This turns out to be a general pattern, which we study next.
\end{exmp}

\begin{say}[Holonomy groups and representations of classical groups]
\label{faithful.reps}
Here we  study vector bundles $E$
over a  projective variety $X$
 whose
holonomy group is contained in an  irreducible  representation
of a product 
$$
\rho:G=G_0\times \prod_{i=1}^m G_i\onto H\supset H_x(E).
$$
where  $G_0$ is a subgroups of the scalars in $GL(E_x)$
and for $i\geq 1$, 
 $G_i$  is one of the classical groups $SL_{n_i}, Sp_{n_i},SO_{n_i}$.
Thus $\rho $ can be obtained from the basic representations
of the $G_i$ by a tensor product of Schur functors ${\mathbb S}_i$.

The easy case is when $\rho$ is an isomorphism.
In this case $E$ corresponds to a principal
$G$-bundle and the basic representation
of each $G_i$ gives  a vector bundle $F_i$ of rank $n_i$
with structure group $G_i$. Here $L\cong F_0$ is a line bundle. 
Thus we obtain that
$$
E\cong L\otimes\bigotimes_{i\geq 1} {\mathbb S}_i(F_i).
$$

The situation is more complicated if $\ker\rho\neq 1$.
To $E$ we can associate a principal $H$-bundle and
the obstruction to lift it to a principal $G$-bundle
lies in $H^2_{et}(X, \ker \rho)$.
If $\ker\rho$ is not connected, then this is never zero.

We can improve the situation by replacing
the groups  
$$
SL_{n_i}, Sp_{n_i},SO_{n_i}\qtq{by} 
GL_{n_i}, GSp_{n_i},GSO_{n_i}
$$
 and extending
$G_0$ to all scalars $\c^*$. Let us denote these groups
by $G_i^*$.
Set $G^*:=\prod_{i\geq 0} G^*_i$ and extend $\rho$
to $\rho^*=G^*\to \c^*\cdot H$. 

Since $G_0$ maps isomorphically onto the scalars, we
see that
$$
\ker \rho^* \cong \prod_{i\geq 1} Z(G^*_i)\cong (\c^*)^m.
$$
The obstruction to lift a principal $\c^*\cdot H$-bundle
to a principal $G^*$-bundle is now
in the Brauer group $Br(X):=H^2_{et}(X, \o_X^*)$.

Therefore, if the Brauer group is zero, then we can lift
our principal $H$-bundle to a principal $G^*$-bundle.
Thus, as before,  we obtain the following:
\end{say}

\begin{prop}\label{Gl-reps.brauer}  Let $X$ be a smooth projective variety
such that $H^2(X,\o_X)=0$ and $H^3(X,\z)$ is torsion free.

Let $E$ be a stable  vector bundle of rank $N$ on $X$ whose
holonomy group is contained in the image of an  irreducible  representation
$$
\rho:G=\prod_{i=1}^m G_i\to GL_N
$$
given by Schur functors  ${\mathbb S}_i$,
where each  $G_i$ is one of the groups $GL_{n_i}, GSp_{n_i},GSO_{n_i}$.
Then there are vector bundles $F_i$ of rank $n_i$ with structure group $G_i$,
 and a line bundle
$L$ such that
$$
E\cong L\otimes\bigotimes_i {\mathbb S}_i(F_i).\qed
$$
\end{prop}

Next we illustrate 
the principle (\ref{hsm.princ}) by studying the
possible holonomy groups and the corresponding
geometric structures for bundles of small rank.

\begin{say}[Holonomy groups of low rank bundles]\label{low.rk.class}
Let $X$ be a smooth projective variety and $E$ a vector bundle
which is stable with respect to an ample divisor $H$.
Assume that the holonomy group
$H_x(E)$ is connected.
\medskip

{\it Rank 2 bundles.} Here $H_x(E)$ is $SL_2$ or $GL_2$.
The first case corresponds to $\det E\cong \o_X$ and
the second to the case when $\det E$ is a line bundle which is not
torsion in  $\pic X$.

The case when $\det E$  is 
torsion in  $\pic X$ would give
nonconnected holonomy.
\medskip

{\it Rank 3 bundles.}
The general case is when $H_x(E)$ is $SL_3$ or $GL_3$.
We can also have $SO_3$,  when $E\cong E^*$
or $GSO_3$ when $E\cong E^*\otimes L$ for some line bundle $L$.

The isomorphism  $SO_3\cong PSL_2$ 
was studied in (\ref{s03.exmp}).
\medskip

{\it Rank 4 bundles.} 
The  cases when $H_x(E)$ is $SL_4, GL_4$
or  $SO_4, GSO_4$ are as before.
We can also have $Sp_4$ or $GSp_4$ holonomy,
corresponding to the existence of a skew symmetric
pairing $E\times E^*\to \o_X$ or 
 $E\times E^*\to  L$ for some line bundle $L$.

There are 2 more interesting cases
when $H_x(E)=SL_2$ or $H_x(E)=GL_2/\mu_3$
with the 3rd symmetric power representation.

Assume that  $H_x(E)=SL_2$ with the 3rd symmetric power representation.
Then $H_x(E)\subset Sp_4$. Furthermore, by (\ref{faithful.reps}),
there is a 
rank 2 vector bundle $F$ such that $E\cong S^3F$.

Finally the last case is when $H_x(E)=GL_2/\mu_3$.
This can be treated as in (\ref{Gl-reps.brauer}).
\medskip

{\it Rank 5 bundles.} 
The  cases when $H_x(E)$ is $SL_5, GL_5$
or  $SO_5, GSO_5$ are as before.

There are 2 other cases
when $H_x(E)=PSL_2$ or $H_x(E)=GL_2/\mu_4$
with the 4th symmetric power representation.
We see that $PSL_2\subset SO_5$ and
and $GL_2/\mu_4\subset GSO_5$, so these have
orthogonal structures.
A more detailed study is given in (\ref{Gl-reps.brauer}).
\medskip

{\it Rank 6 bundles.} 
The  cases when $H_x(E)$ is $SL_6, GL_6$,
  $SO_6, GSO_6$ or  $Sp_6$, $GSp_6$ are as before.

The cases when
$H_x(E)=SL_2$ or $H_x(E)=GL_2/\mu_5$
with the 5th symmetric power representation
or
$H_x(E)=SL_3$ or $H_x(E)=GL_3/\mu_2$
with the 2nd symmetric power representation
work as the rank 4 cases.

By luck, the cases when
$H_x(E)=SL_4/\mu_2$ or $H_x(E)=GL_4/\mu_2$
with the 2nd exterior power representation
are contained in $SO_6$ (resp.\ $GSO_6$).

The last case is when $H_x(E)=SL_2\times SL_3$,
or $H_x(E)=GL_2\times GL_3/\c^*$
with the tensor product of the standard
representations.
In the first case, by
(\ref{faithful.reps}),  we get that $E\cong F_2\otimes F_3$
with $H_x(F_2)=SL_2$ and $H_x(F_3)=SL_3$
while in the  second case
we again have a Brauer obstruction to deal with
(\ref{Gl-reps.brauer}).

Every other connected reductive subgroup of
$GL_6$ is contained in one of the above.
\medskip

{\it Rank 7 bundles.} 
The  cases when $H_x(E)$ is $SL_7, GL_7$ or
  $SO_7, GSO_7$  are as before.

The cases when
$H_x(E)=PSL_2$ or $H_x(E)=GL_2/\mu_6$
with the 6rd symmetric power representation
are examined in  (\ref{Gl-reps.brauer}).

The first exceptional case also appears, namely
we can have monodromy group $G_2\subset SO_7$
or $\c^*\cdot G_2\subset GSO_7$. We can not say anything useful about
it beside noting that the monodromy is a subgroup of $\c^*\cdot G_2$
iff    $\wedge^3 E$ has a line bundle direct summand.
This is a  consequence of the
corresponding characterization of 
$G_2\subset SL_7$ as the subgroup that fixes a general 
skew symmetric trilinear form.

Indeed, one checks that in the $SL_7, GL_7, SO_7, GSO_7$ 
cases there is no 1--dimensional
invariant subspace in $\wedge^3 E$.
In the $PSL_2$ or $GL_2/\mu_6$ cases there is such an
invariant subspace, but $PSL_2\subset G_2$ and  
$GL_2/\mu_6\subset \c^*\cdot G_2$.
\medskip

{\it Rank 8 bundles.} Here we get the first case of a bundle
$E$ with $c_1(E)=0$ where the pattern of (\ref{low.rank.prop}) no longer holds.

This is when the holonomy group is the tensor product of the
standard representations of $SL_2$ and $SL_4$.
Thus $H_x(E)\cong (SL_2\times SL_4)/\mu_2$ and the
Brauer obstruction is inevitable.
\medskip

{\it Rank $\leq 16$ bundles.}
By now it should be clear that one can continue in this manner
for low ranks, and either
direct constructions or the method of (\ref{Gl-reps.brauer}) apply.

For rank 16 we run into the first case where the holonomy can be
a spinor group, here $H_x(E)\cong Spin_5$. Probably it is again 
a Brauer--type
obstruction, whose vanishing ensures that 
$E$ is one of the half spin subbundles of the Clifford algebra of a
rank 5 bundle $F$ with orthogonal structure.
\end{say}

\section{Tangent bundles}

It may be especially interesting to consider the
holonomy group of the tangent bundle $T_X$ of a
smooth projective variety $X$.
There are only a few cases when
we can compute the algebraic holonomy group.

\begin{say}[Calabi--Yau varieties] Let $X$ be a 
simply connected smooth projective variety $X$
such that $K_X=0$ which is not a direct product.

The differential geometric holonomy group is 
either $SU_n(\c)$ or $U_n(\h)$. As observed in \cite{beauville}, the
tensor powers  of $T_X$ decompose according to
the representation theory of $SU_n$ (resp.\ or $U_n(\h)$), 
thus by (\ref{main.thm.rem}.3) we conclude that the algebraic holonomy
is $SL_n(\c)$ (resp.\ $Sp_{2n}(\c)$).
This proves (\ref{hol=hol.CY}).
\end{say}

\begin{say}[Homogeneous spaces] Let $X=G/P$ be a
smooth, projective homogeneous space under a reductive group $G$.
Let $\rho:P\to GL(T_xX)$ denote the stabilizer representation.
The stabilizer representation vanishes on the unipotent radical
$U\subset P$ and so we can view $\rho$ as a 
representation of the reductive Levi subgroup
$\rho:P/U\to GL(T_xX)$.

The tangent bundle $T_X$ is indecomposable iff $\rho$ is irreducible.
By \cite{ram-homog, ume, kobayashi}, in this case  $T_X$ is stable 
and  tensor powers  of $T_X$ decompose according to
the representation theory of $P/U$. 
Thus by (\ref{main.thm.rem}.3) we conclude that the algebraic holonomy
group is  $\rho(P/U)$. This proves (\ref{G/P.prop}).

There are very few examples of Fano varieties
 whose holonomy group  we can compute.

For instance, let $S$ be a Del Pezzo surface which is
obtained from $\p^2$ by blowing up at least 3 points.
It is easy to see that $T_S$ is stable,
hence by (\ref{hol.simpconn.prop}) the holonomy group is $GL_2$.

\begin{ques} Let $X^n$ be a  smooth projective variety 
 with Picard number 1 and $-K_X$  ample.
Assume that the automorphism group of $X$ is finite
and $T_X$ is stable. Is the 
algebraic holonomy group  $GL_n$?
\end{ques}

\end{say}

\begin{say}[Varieties with ample canonical class]
Let $X$ be a  smooth projective variety $X$
such that $K_X$ is ample. 

By the Akizuki--Nakano vanishing theorem
(cf.\ \cite[p.155]{g-h}),
$\wedge^i\Omega_X$ contains no ample line bundle for $i<\dim X$.
In particular, $\Omega_X$ does not contain any subsheaf of
rank $<\dim X$ whose determinant is ample and so 
 it is stable with respect to the ample divisor
$K_X$. (This also follows from the much stronger
result of Aubin and Yau about the existence of  a K\"ahler--Einstein metric.)

Furthermore, this also implies that 
$\wedge^i\Omega_X$ has no  line bundle direct summands for $i<\dim X$.
Thus we conclude:

\begin{prop}\label{gen.type.AN.hol} Let $X$ be a  smooth projective variety $X$
such that $K_X$ is ample and let $H_x\subset GL(T_xX)$ denote
the holonomy group of the tangent bundle $T_X$.

Then $\wedge^iT_xX$ has no  1--dimensional $H_x$-invariant
subspaces for $i<\dim X$. In particular, $H_x$ acts
irreducibly on $T_xX$.
\end{prop}

Thus it is natural to study subgroups $H\subset GL(E)$
such that $\wedge^iE$ has no  1--dimensional $H$-invariant
subspace for $i<\dim E$. This is a very restrictive condition,
but we have not been able to classify all such
representations.
In any case, at the moment we do not even know 
the answer to the following:

\begin{ques} Is there  a  simply connected, smooth projective variety $X$
 with Picard number 1 and $K_X$  ample, whose
algebraic holonomy group is different from $GL_n$?
\end{ques}

There are smooth projective varieties 
 with Picard number 1 which are quotients of a direct product,
and these have smaller holonomy group.
\end{say}

 \begin{ack}  We thank F.\ Bogomolov, R.\ Guralnick, R.\ Lazarsfeld, 
R.\ Livn\'e, P.\ Sarnak 
and C.S \ Seshadri
   for useful comments, references  and suggestions. 
We are grateful for   many discussions with  M.\ Larsen concerning
representation theory.
 Partial financial support for
   JK was provided by the NSF under grant number DMS-0500198.
\end{ack}

\bibliography{refs}

\bigskip

\noindent Chennai Math.\ Inst.\ SIPCOT IT Park, Siruseri-603103, India

\begin{verbatim}balaji@cmi.ac.in\end{verbatim}

\bigskip

\noindent Princeton University, Princeton NJ 08544-1000, USA

\begin{verbatim}kollar@math.princeton.edu\end{verbatim}

\end{document}